\documentclass{gtpart}     
\usepackage{pinlabel}  
\usepackage{graphicx}  
\usepackage[all]{xy}
\usepackage{amscd}
\usepackage{enumerate}
\usepackage{fullpage}
\title{Cohomology with twisted coefficients of the classifying space of a fusion system}

\author{R\'emi Molinier}
\givenname{R\'emi}
\surname{Molinier}
\address{Universit\'e Paris 13, Sorbonne Paris Cit\'e, LAGA, UMR 7539 du CNRS, 99, Av. J.-B. Cl\'ement, 93430 Villetaneuse, France.}
\email{molinier@math.univ-paris13.fr}
\urladdr{https://www.math.univ-paris13.fr/~molinier/index.html}


\arxivreference{}
\arxivpassword{}

\volumenumber{}
\issuenumber{}
\publicationyear{}
\papernumber{}
\startpage{}
\endpage{}
\doi{}
\MR{}
\Zbl{}
\received{}
\revised{}
\accepted{}
\published{}
\publishedonline{}
\proposed{}
\seconded{}
\corresponding{}
\editor{}
\version{}

\newtheorem{thm}{Theorem}[section]    
\newtheorem{lem}[thm]{Lemma}          
\newtheorem{prop}[thm]{Proposition}
\newtheorem{cor}[thm]{Corollary} 
\newtheorem{thmx}{Theorem}

\theoremstyle{definition}
\newtheorem{defi}[thm]{Definition}
\newtheorem{conj}[thm]{Conjecture}

\newtheorem{rem}[thm]{Remark}             

\newtheorem*{nota}{Notation}
\newcommand{\F}{\mathbb{F}}
\newcommand{\Zp}{\mathbb{Z}_{(p)}}
\newcommand{\Mod}{\text{-Mod}}

\newcommand{\Fc}{\mathcal{F}}
\newcommand{\Li}{\mathcal{L}}
\newcommand{\T}{\mathcal{T}}
\newcommand{\Bc}{\mathcal{B}}
\newcommand{\A}{\mathcal{A}}
\newcommand{\B}{\mathcal{B}}
\newcommand{\Dc}{\mathcal{D}}
\newcommand{\Ss}{\mathcal{S}}
\newcommand{\Ch}{\text{Ch}}

\newcommand{\pcompl}[1]{#1^{\wedge}_p}
\newcommand{\limproj}[1]{\lim\limits_{\substack{\longleftarrow \\ #1}}}

\newcommand{\Ob}{\text{Ob}}
\newcommand{\Hom}{\text{Hom}}
\newcommand{\Aut}{\text{Aut}}

\newcommand{\Inj}{\text{Inj}}

\newcommand{\Mor}{\text{Mor}}
\newcommand{\End}{\text{End}}

\newcommand{\incl}{\text{incl}}
\newcommand{\Res}{\text{Res}}
\newcommand{\tr}{\text{tr}}

\newcommand{\Id}{\text{Id}}
\newcommand{\im}{\text{Im}\,}

\newcommand{\Syl}{\text{Syl}}

\begin{document}

\begin{abstract}    
We study the cohomology with twisted coefficients of the geometric realization of a linking system associated to a saturated fusion system $\Fc$.
More precisely, we extend a result due to Broto, Levi and Oliver to twisted coefficients.
We generalize the notion of $\Fc$-stable elements to $\Fc^c$-stable elements in a setting of cohomology with twisted coefficients by an action of the fundamental group.
We then study the problem of inducing an idempotent from an $\Fc$-characteristic $(S,S)$-biset 
and we show that, if the coefficient module is nilpotent, then the cohomology of the geometric realization of a linking system
can be computed by $\Fc^c$-stable elements. 
As a corollary, we show that for any coefficient module, 
the cohomology of the classifying space of a $p$-local finite group
can be computed by these $\Fc^c$-stable elements.  
\end{abstract}

\maketitle


The notion of saturated fusion systems was introduced by Puig in the 90s in a context of modular representation theory. In topology, saturated fusion systems are used in the study of $p$-completed classifying spaces of groups.
A \emph{$p$-local finite group} is a triple $(S,\Fc,\Li)$ where $S$ is a $p$-group, $\Fc$ a saturated fusion system over $S$ and $\Li$ an associated centric linking system.
For a $p$-local finite group $(S,\Fc,\Li)$, the $p$-completed nerve of $\Li$, $\pcompl{|\Li|}$, is called its \emph{classifying space}. 
The theory of $p$-local finite group has been studied in detail by Broto, Levi, Oliver and others 
(see \cite{BLO2}, \cite{OV1},\cite{5a1} and \cite{5a2}). 
The linking system and its geometric realization, even without $p$-completion, play here a fundamental and central role. 
In fact, for a given saturated fusion system, the existence and uniqueness of a linking system associated to a saturated fusion system were shown more recently by Chermak \cite{Ch} (using the theory of partial groups). The proof of this important conjecture highlights the deep link between group theory and homotopy theory and the importance of linking systems to study $p$-local structures. We refer to Aschbacher, Kessar and Oliver \cite{AKO} for more details about fusion systems in general.

A well-known result of Cartan and Eilenberg (see \cite{CE} Chap XII, Theorem 10.1) 
expresses the cohomology with mod $p$ coefficients of a finite group as the submodule of \emph{stable} elements in the cohomology of a Sylow $p$-subgroup.
This submodule of stable elements corresponds to the inverse limit over the fusion system of the group cohomology functor.
One important result in the theory of $p$-local finite groups is an extension of this theorem to any $p$-local finite group 
which tells us that the cohomology of the geometric realization of a linking system with trivial coefficients can be computed by $\Fc$-stable elements.

\begin{thmx}[\cite{BLO2}, Theorem B]\label{trivial}
Let $(S,\Fc,\Li)$ be a $p$-local finite group.
The inclusion of $BS$ in $|\Li|$ induces a natural isomorphism
\[
 \xymatrix{ H^*(\pcompl{|\Li|},\F_p)\cong H^*(|\Li|,\F_p)\ar[r]^-\cong & H^*(\Fc^c,\F_p):=\limproj{\Fc^c} H^*(-,\F_p).}
\]
\end{thmx}

Here $\Fc^c$ denotes the full subcategory of $\Fc$ where the objects are the subgroups of $S$ which are $\Fc$-centric 
(the analog of $p$-centric in the group case).
From that, the above result can be extended to cohomology with coefficients in any $\Zp$-module $A$ (a proof is given in \cite{5a2} Lemma 6.12).

The Theorem of Cartan and Eilenberg actually applies with other choices of coefficients and especially with twisted coefficients. 
One question asked by Oliver in his book with Aschbacher and  Kessar \cite{AKO} is the understanding of the cohomology of $|\Li|$ with twisted coefficients. 
Indeed, this cohomology appears for example in the study of extensions of $p$-local finite groups \cite{5a2} and it can give more information about the link between the fusion system and the homotopy type of $|\Li|$.

In this paper, cohomology with twisted coefficients means cohomology with locally constant coefficients.
In other words, given a space $X$, the cohomology of $X$ with twisted coefficients corresponds to the cohomology of $X$ with coefficients in a $\Zp[\pi_1(X)]$-module.
Let $M$ be a $\Zp[\pi_1(X)]$-module, and suppose in addition that $X$ has a universal covering space $\widetilde{X}$, the cohomology of $X$ with twisted coefficients in $M$
is the cohomology of the chain complex
\[C^*(X;M)=\Hom_{\Z[\pi_1(X)]}(S_*(\widetilde{X}),M),\]
where $S_*(\widetilde{X})$ is the usual singular chain complex of $\widetilde{X}$.
 
In this paper we extend Theorem \ref{trivial} to twisted coefficients. 
Levi and Ragnarsson allready consider this problem in \cite{LR}. As they pointed out (\cite{LR}, Proposition 3.1), it is not possible to achieve a version of Theorem \ref{trivial} for twisted coefficients in full generality and some restrictions on $\Fc$ or on the coefficients are necessary. 

\medskip
{\bf Theorem \ref{main}}\qua
{\sl  Let $(S,\Fc,\Li)$ be a $p$-local finite group.
 If $M$ is an abelian $p$-group with a nilpotent action of $\pi_1(|\Li|)$, then the inclusion of $BS$ in $|\Li|$ induces a natural isomorphism,
\[ H^*(|\Li|,M)\cong H^*(\Fc^c,M).\] }
Here, for $p$ a prime, a $p$-group is a finite group of $p$-power order.
An important application of the above result is the complete description of the cohomology of the classifying space of a fusion system in terms of $\Fc^c$-stable elements.

\medskip
{\bf Corollary \ref{main2}}\qua
{\sl Let $(S,\Fc,\Li)$ be a $p$-local finite group.
If $M$ is an abelian $p$-group with an action of $\pi_1(\pcompl{|\Li|})$, then there is a natural isomorphism
\[
 H^*(\pcompl{|\Li|},M)\cong H^*(\Fc^c,M).
\]}

One crucial tool in the proof of Theorem \ref{main}, as it was in the proof of Theorem \ref{trivial}, is to construct an idempotent of $H^*(S,M)$ from an $\Fc$-characteristic $(S,S)$-biset which generalizes the construction of Broto, Levi and Oliver in \cite{BLO2}. The action on the coefficients makes the construction less straightforward and we need here to work with a finite module to guarantee that this idempotent exists. 

\vspace{10pt}

\textbf{Organization.}
We start with a review on $p$-local finite groups in Section \ref{2}.
In Section \ref{3}, we introduce our cohomology functor and 
we define properly the notion of $\Fc^c$-stable elements.
Section \ref{4} contains the construction of an idempotent for an $\Fc$-characteristic biset. As this section is rather long, it is divided into several subsections, each of which deals with a different aspect or property of the idempotent.
We prove Theorem \ref{main} in the Section \ref{5}.
Finally, in Section \ref{6}, we give a result about the cohomology with twisted coefficients of $p$-good spaces and we apply it to get Corollary \ref{main2}.

\textbf{Acknowledgments.}
I first would like to thank Bob Oliver, my PhD supervisor, for his support during my thesis and this work.
I also would like to thank Radu Stancu for all the time we spent together on this problem. I am grateful to the Center for Symmetry and Deformation for its hospitality and Jesper Grodal and Sune Precht Reeh for many fruitful conversations.

Finally, I would like to thank the referee, Alex Gonz\'ales and Emilie Devijver for their careful reading of the paper and their useful comments.

This work was supported by the Universit\'e Paris 13.

\section{Background on \texorpdfstring{$p$}{Lg}-local finite groups}\label{2}

We give here a very short introduction to $p$-local finite groups. 
 We refer the reader interested in more details to Aschbacher, Kessar and Oliver \cite{AKO}.

Roughly speaking, fusion systems encode the conjugation data of a finite group with respect to a choice of a Sylow $p$-subgroup.
For $G$ a finite group and $g\in G$, we will denote by $c_g$ the homomorphism $x\in G\mapsto gxg^{-1}\in G$. Given subgroups $H,K\leq G$, we shall denote by $\Hom_G(H,K)$  the set of all group homomorphisms $c_g$ for $g\in G$ such that $c_g(H)\leq K$.

\begin{defi}\label{defF}
 Let $S$ be a finite $p$-group.
 A \emph{fusion system} over $S$ is a small category $\Fc$, where $\Ob(\Fc)$ is the set of all subgroups of $S$ and which satisfies the following two
 properties for all $P,Q\leq S$:
 \begin{enumerate}[(a)]
  \item $\Hom_S(P,Q)\subseteq \Mor_\Fc(P,Q)\subseteq \Inj(P,Q)$;
  \item each $\varphi\in\Mor_\Fc(P,Q)$ is the composite of an $\Fc$-isomorphism followed by an inclusion.
 \end{enumerate}
A fusion system is \emph{saturated} if it satisfies two more technical conditions called the saturation axioms (we refer the reader to \cite{AKO}, Definition I.2.1 for a proper definition).
\end{defi}

The composition in a fusion system is given by composition of homomorphisms. We usually write $\Hom_\Fc(P,Q)=\Mor_\Fc(P,Q)$ to emphasize the fact that the 
morphims in $\Fc$ are actual group homomorphisms.

The typical example of saturated fusion system is the fusion system $\Fc_S(G)$ of a finite group $G$ over $S\in\Syl_p(G)$. 

\begin{defi}\label{pcentr}
Let $\Fc$ be a saturated fusion system over a $p$-group $S$.
A subgroup $P\leq S$ is \emph{$\Fc$-centric} if $C_S(Q)=Z(Q)$ for every $Q\in P^\Fc$.
We will denote by $\Fc^c$ the full subcategory of $\Fc$ with set of objects all the $\Fc$-centric subgroups of $S$.
\end{defi}

If $\Fc$ is the saturated fusion system associated to a finite group $G$ with $S$ as Sylow $p$-subgroup, then a subgroup $P\leq S$ is $\Fc$-centric if and only if $P$ is $p$-centric, i.e. $Z(P)$ is a Sylow $p$-subgroup of $C_G(P)$.
Before defining the notion of centric linking system let us first recall a well-known result about saturated fusion system.

\begin{thm}[Alperin's Fusion Theorem]\label{AFT}
Let $\Fc$ be a saturated fusion system over a $p$-group $S$.
Then, every morphism is a composite of restrictions of automorphisms of $\Fc$-centric subgroups.  
\end{thm}

In other words, a saturated fusion system $\Fc$ is generated by $\Fc^c$. 
In fact, Alperin's Fusion Theorem is more precise and says that we just need automorphisms of $S$ and $\Fc$-essential subgroups of $S$.
For more details, we refer to Section I.3 of Aschbacher, Kessar and Oliver \cite{AKO}.

For $\Fc$ a fusion system over a $p$-group, $\T_S^c(S)$ will denote the usual transporter category of $S$ with set of objects $\Ob(\Fc^c)$.
\begin{defi}\label{linkdef}
 Let $\Fc$ be a fusion system over a $p$-group $S$.
 A \emph{centric linking system} associated to $\Fc$ is a finite category $\Li$ together with a pair of functors
 \[
  \xymatrix{\T_S^c(S)\ar[r]^-{\delta} &\Li \ar[r]^-{\pi}& \Fc} 
 \]
satisfying the following conditions:
\begin{itemize}
 \item[(A)] $\delta$ is the identity on objects, and $\pi$ is the inclusion on objects. For each $P,Q\in\Ob(\Li)$ such that $P$ is 
       fully centralized in $\Fc$, $C_S(P)$ acts freely on $\Mor_\Li(P,Q)$ via $\delta_{P,P}$ and right composition, and 
       \[
        \xymatrix{\pi_{P,Q}:\Mor_\Li(P,Q)\ar[r] & \Hom_\Fc(P,Q)}
       \]
       is the orbit map for this action.
 \item[(B)] For each $P,Q\in\Ob(\Li)$ and each $g\in T_S(P,Q)$, the application $\pi_{P,Q}$ sends $\delta_{P,Q}(g)\in\Mor_\Li(P,Q)$ to $c_g\in\Hom_\Fc(P,Q)$.
 \item[(C)] For each $P,Q\in\Ob(\Li)$, all $\psi\in\Mor_\Li(P,Q)$ and all $g\in P$, the diagram
 \[
  \xymatrix{ P\ar[d]_ {\delta_P(g)} \ar[r]^\psi & Q \ar[d]^{\delta_Q(\psi(g))} \\
             P \ar[r]^\psi & Q}
 \]
 commutes in $\Li$.
\end{itemize}
 A \emph{$p$-local finite group} is a triple $(S,\Fc,\Li)$ where $S$ is a $p$-group, $\Fc$ a saturated fusion system over $S$,
 and $\Li$ is a linking system associated to $\Fc$.
The \emph{classifying space} of $(S,\Fc,\Li)$ is then given by $\pcompl{|\Li|}$. 
\end{defi}


\section{Cohomology and \texorpdfstring{$\Fc^c$}{Lg}-stable elements}\label{3} 

In this section we introduce our cohomology functor with twisted coefficients defined on $\Fc^c$ and we define the notion of $\Fc^c$-stable elements. We refer the reader to \cite{We} for all the necessary results on homological algebra as well as the classicial notion of $\delta$-functor. We also refer the reader to \cite{We}, \cite{Br}, \cite{CE} or \cite{AM} for details on group cohomology.
\vspace{10pt}

Before introducing the notion of $\Fc^c$-stable elements, we need to understand the action of $S$ on a $\Zp[\pi_1(|\Li|)]$-module.
For each pair of $\Fc$-centric subgroups $P\leq Q$, set
\[\iota_P^Q=\delta_P^Q(1).\] 
We denote by
\[\pi_\Li=\pi_1(|\Li|,S)
\]
the fundamental group of the geometric realization $|\Li|$ with base point at the vertex $S$. 
For $G$ a discrete group, let $\Bc(G)$ be the category with 
a unique object, and morphism set equals to $G$ (hence, $|\Bc(G)|=BG$).
Consider the functor 
\[\xymatrix{\omega:\Li\ar[r] & \Bc(\pi_\Li)}\]
which maps each object to the unique object in the target and sends each morphism $\varphi\in\Mor_\Li(P,Q)$ to the class of the loop 
$\iota_Q.\varphi. \overline{\iota_P}$ where $\overline{\iota_P}$ is the edge $\iota_P$ followed in the opposite direction.
In particular, every $\Zp[\pi_\Li]$-module is naturally a $\Zp[S]$-module where the action is given by the following composition:
\[\xymatrix{\Bc(S)=\Bc(\Mor_{\T^c_S(S)}(S,S))\ar[r]^-{\delta_S} & \Li \ar[r]^-\omega & \Bc(\pi_\Li)}.\] 
\vspace{10pt}

 

Let $(S,\Fc,\Li)$ be a $p$-local finite group and let $M$ be a $\Zp[\pi_\Li]$-module. As we work with an action of $\pi_1(|\Li|)$, 
we can define a functor on $\Li$ using the bi-functoriality of group cohomology. 
Recall that group cohomology defines a contravariant functor
\[\xymatrix{H^*(-,-):\Dc\ar[r] & \Zp\Mod}\] 
where $\Dc$  is the category of pairs $(G,M)$ with $G$ a group and $M$ a $\Zp[G]$-module. A morphism in $\Dc$ from $(G,M)$ to $(H,N)$ is a pair, $(\varphi,\rho)$ where $\xymatrix{\varphi:G\ar[r]& H}$ is a group homomorphism and $\xymatrix{\rho: N\ar[r] & M}$ is a linear map such that, for every $n\in N$ and 
every $g\in G$, $g\rho(n)=\rho\left(\varphi(g)n\right)$.

Given $\varphi\in\Mor(\Li)$, we have $(\pi(\varphi),\omega(\Li))\in\Mor(\Dc)$ and we define our cohomology functor as the following.
\[
 \xymatrix@R=1mm @C=1cm{**[r]H^*(-,M):& \Li \ar[r] & **[r]\Zp\Mod \\
                    & P\in\Ob(\Li) \ar@{|->}[r] & **[r]H^*(P,M)\\
                    & \varphi\in\Mor_\Li(P,Q)\ar@{|->}[r] & **[r] H^*(\varphi,M)=\varphi^*:=H^*(\pi(\varphi),\omega(\varphi)^{-1}).} 
\]

For $P,Q$ two subgroups of $S$ and $\varphi\in\Mor_\Li(P,Q)$, $H^*(\varphi,M)$ can also be defined on the chain level as follows:
 \[
 \xymatrix@R=1mm{\Hom_{\Zp[Q]}\left(R_\bullet,M\right) \ar[r] & \Hom_{\Zp[P]}\left(R_\bullet,M\right)\\
              f  \ar@{|->}[r] &   \left( \omega(\varphi)^{-1}\circ f\circ\pi(\varphi)_*\right)}
 \]
where $(R_\bullet)$ is a projective resolution of the trivial $\Zp[S]$-module $\Zp$.
Finally, it can also be defined as the morphism between the two derived functors of $(-)^Q$ and $(-)^P$ induced by 
\[
\xymatrix{ x\in M^Q \ar@{|->}[r] & \omega(\varphi)^{-1}x \in M^P.}
\]
By construction, it defines a morphism of $\delta$-functors.

\begin{prop}\label{phimorhdelta}
 Let $(S,\Fc,\Li)$ be a $p$-local finite group.
 If $P,Q\leq S$ are $\Fc$-centric and $\varphi\in\Mor_\Li(P,Q)$, then 
 \[
  \xymatrix{H^*(\varphi,-): H^*(Q,-)\ar[r] & H^*(P,-)}
 \] 
 is a morphism of $\delta$-functors from $\left(H^*(Q,-),\delta_{H^*(Q,-)}\right)$ to $\left(H^*(P,-),\delta_{H^*(P,-)}\right)$. 
\end{prop}


By construction, this functor naturally extends the group cohomology functor defined on $\T_S^c(S)$.
\[
\xymatrix{ \T_S^c(S)\ar[rr]^-{H^*(-,M)}\ar[rd]_\delta& & \Zp\Mod \\
                 & \Li \ar[ru]_-{H^*(-,M)}& }          
\]
In particular, for every $P\leq S$ and $g\in P$, $H^*(\delta_P(g),M)=c_g^*$.
Moreover, it also factors through $\Fc^c$ along $\xymatrix{\pi:\Li\ar[r]& \Fc^c}$.

\begin{prop}
Let $\varphi,\beta\in\Mor_\Li(P,Q)$ with $P,Q\in\Li$.
If $\pi(\varphi)=\pi(\beta)$ then $H^*(\varphi,M)=H^*(\beta,M)$.
\end{prop}

\begin{proof}
If $\pi(\varphi)=\pi(\beta)$, then there exists $u\in Z(P)$ such that $\varphi=\beta\circ \delta_P(u)$ and thus 
\[H^*(\varphi,M)=H^*(\delta_P(u),M)\circ H^*(\beta,M).\]
However $H^*(\delta_P(u),M)= H^*(\pi(\delta_P(u)),\omega(u)^{-1})=H^*(c_u,\omega(u)^{-1})=c_u^*$ 
is the automorphism of $H^*(P,M)$ induced by the conjugation by $u$, and, as $u\in Z(P)\leq P$, 
this is the identity.
\end{proof}
In particular, if $\pi(\varphi)=\incl^Q_P$, then $H^*(\varphi,M)=H^*(\iota^Q_P,M)=H^*(\incl^Q_P,\Id_M)=\Res^Q_P$. Hence, $H^*(-,M)$ factors naturally through $\Fc^c$ along $\pi$. 
For $M$ a $\Zp[\pi_\Li]$-module, $P,Q\leq S$ two $\Fc$-centric subgroups and $\varphi\in\Hom_\Fc(P,Q)$ we write $\varphi^*:=H^*(\psi,M)$ 
where $\psi\in\Mor_\Li(P,Q)$ is such that $\pi(\psi)=\varphi$.

\begin{defi}
 An element $x\in H^*(S,M)$ is called \emph{$\Fc$-centric stable}, or just \emph{$\Fc^c$-stable}, if for all $P\in\Ob(\Fc^c)$ and all $\varphi\in\Hom_\Fc(P,S)$,
 \[\varphi^*(x)=\Res_P^S(x).\]
We denote by $H^*(\Fc^c,M)\subseteq H^*(S,M)$ the submodule of all $\Fc^c$-stable elements.
\end{defi}

This submodule of $\Fc^c$-stable elements corresponds to the inverse limit of $H^*(-,M)$ on the category $\Fc^c$,
\[H^*(\Fc^c,M)\cong\limproj{\Fc^c} H^*(-,M).\]
Moreover, if $M$ is a $\Zp$-module with a trivial action of $\pi_\Li$, then, by Alperin's Fusion Theorem (Theorem \ref{AFT}),
$H^*(\Fc^c,M)\cong H^*(\Fc,M)$ and the notion of $\Fc^c$-stable elements naturally extends the notion of $\Fc$-stable elements.

In general, we cannot expect to define a cohomology functor on all $\Fc$. For example, every morphism of $\Li$ induces the identity on the trivial subgroup $\{e\}$. Thus, if the cohomology functor was defined on all $\Fc$,  every morphism in $\Li$ should act trivially on $M=H^0(\{e\},M)$, which is absurd.
Instead we consider the following construction.

Let $(S,\Fc,\Li)$ be a $p$-local finite group.
Let $Q\leq S$ be a $\Fc$-centric subgroup of $S$, $P_0\leq Q$, $\psi\in\Aut_\Li(Q)$ and denote $P_1=\pi(\psi)(P_0)\leq Q$. 
For $M$ a $\Zp[\pi_\Li]$-module, even if $P_0$ is not $\Fc$-centric, we can consider the morphism $\xymatrix{H^*(P_1,M)\ar[r]& H^*(P_0,M)}$ 
given by $H^*(\pi(\psi)|^{P_1}_{P_0},\omega(\psi))$. This can also be defined on the chain level by,
 \[
 \xymatrix@R=1mm{\Hom_{\Zp[P_1]}\left(R_\bullet,M\right) \ar[r] & \Hom_{\Zp[P_0]}\left(R_\bullet,M\right)\\
              f  \ar@{|->}[r] &   \left( \omega(\psi)^{-1}\circ f\circ(\pi(\psi)|^{P_1}_{P_0})_*\right)}
 \]
where $(R_\bullet)$ is a projective resolution of the trivial $\Zp[S]$-module $\Zp$.

This is well-defined because, by Definition \ref{linkdef} $(C)$, 
\[\xymatrix @R=1mm { M \ar[r] & M \\
x\ar@{|->}[r] & \omega(\psi)^{-1}x}\] 
defines a linear map such that, for every $p\in P_0$ and $x\in M$, 
\[
\omega(\psi)^{-1} \omega\left(\delta_Q(p)\right) x=\omega\left(\delta_Q(\psi(p))\right)	\omega(\psi)^{-1}x
\]
and thus $(\pi(\psi)|^{P_1}_{P_0},\omega(\psi))\in\Mor(\Dc)$.
Notice that if $P$ is $\Fc$-centric, then this is just $H^*(\psi,M)$.
Note also that, if the action of $\omega(\psi)$ on $M$ is trivial, this is just the usual morphism induced in cohomology by $\pi(\psi)|^{P_1}_{P_0}$.

Let $P,Q\leq S$ and $\varphi\in\Hom_\Fc(P,Q)$. By Alperin's Fusion Theorem \ref{AFT}, there exist $P=P_0, P_1,\dots,P_r=\varphi(P)$ subgroups of $S$, $Q_1,\dots,Q_r$ $\Fc$-centric subgroups of $S$ and $\psi_i\in\Aut_\Li(Q_i)$ for every $i$ such that $\pi(\psi_i)(P_{i-1})=P_{i}$ and 
\[\varphi=\pi(\psi_r)|^{P_r}_{P_{r-1}}\circ\pi(\psi_{r-1})|^{P_{r-1}}_{P_{r-2}}\circ\dots\circ\pi(\psi_1)|^{P_1}_{P_0}.\]
We then consider the following composite
\[
 \xymatrix@C=3cm{ H^*(P_r,M) \ar[r]^-{H^*(\pi(\psi_1)|^{P_1}_{P_0},\omega(\psi_1)^{-1})} &  \cdots \ar[r]^-{H^*(\pi(\psi_r)|^{P_r}_{P_{r-1}},\omega(\psi_r)^{-1})} & H^*(P,M)} 
\]
composed on the right by $Res^Q_{P_r}$ which gives us a morphism 
\[\xymatrix{H^*(Q,M)\ar[r] & H^*(P,M).}\]
Note that this morphism depends on the choice of the decomposition of $\varphi$ into restrictions of automorphisms of $\Fc$-centric subgroups and
not only on $\varphi$. As an example, we can again look at the trivial subgroup $\{e\}$ in a given fusion system $\Fc$: 
each morphism in $\Fc^c$ restricts to the identity on $\{e\}$, but, if $M$ is not a trivial $\Zp[\pi_\Li]$-module,
not every $\varphi\in\Mor(\Fc^c)$ acts trivially on $M=M^{\{e\}}=H^0(\{e\},M)$.

\begin{rem}\label{phimorhdelta2}
By construction, a morphism $\xymatrix{H^*(Q,M)\ar[r] & H^*(P,M)}$ obtained from $\varphi$ by this process defines a morphism of $\delta$-functors.
\end{rem}

\section{Bisets and idempotents}\label{4}

 An important result in Broto, Levi and Oliver \cite{BLO2}, and a crucial tool in the proof of Theorem \ref{trivial},
 is the existence of an $\Fc$-characteristic $(S,S)$-biset which leads to an idempotent of $H^*(S,\F_p)$ whose image is $H^*(\Fc,\F_p)$.

This section is divided into several subsections. We start by recalling some results about left-free bisets. Next, we describe the interaction of bisets with cohomology with trivial coefficients. The third subsection analyzes the relation of bisets with nontrivial coefficients. After that, we present the construction of an idempotent from an $\Fc$-characteristic biset. As an example, we consider in a fifth subsection the particular case of constrained fusion systems. Finally, the last section gives the link with $\delta$-functors. 
 
\subsection{Background on bisets}\label{4.1}

Let $G,H$ be two finite groups.
Transitive $(G,H)$-bisets (here, $G$ acts on the left and $H$ on the right) are isomorphic to bisets of the form $(G\times H)/K$ for $K$ a subgroup of $G\times H$. 
We can then use the Goursat Lemma to describe all these subgroups.
Here, we are just interested in isomorphism classes of $(G,H)$-bisets where the action of $G$ is free. 
In this setting, the classes of transitive left-free $(G,H)$-bisets are given by pairs $(K,\varphi)$, 
where $K$ is a subgroup of $G$ and $\varphi\in\Hom(K,H)$ is a group homomorphism.

\begin{nota}
For all $(K,\varphi)$, with $K$ a subgroup of $G$ and $\varphi\in\Hom(K,H)$ a group homomorphism, we write
\[\Delta(K,\varphi)=\left\{(k,\varphi(k))\;;\; k\in K\right\}\leq G\times H.\]
 
For a $(G,H)$-pair $(K,\varphi)$, the set $\lbrace K,\varphi\rbrace:=(G\times H)/\Delta(K,\varphi)$ defines a $(G,H)$-biset.
Moreover, its isomorphic class is determined by the conjugacy class of $\Delta(K,\varphi)$ and we denote by $[K,\varphi]$ this class.
\end{nota}

We can also define a category $\Bc$, often called the \textit{Burnside category}, where the objects are the finite groups and, 
for all finite groups $G$ and $H$, $\Bc(G,H)$ is the set of isomorphism classes of $(G,H)$-bisets.
The composition is given by the following construction.

\begin{defi} \label{compdef}
Let $G,H$ and $K$ be finite groups, $\Omega$ a $(G,H)$-biset and $\Lambda$ a $(H,K)$-biset.
We define, 
\[\Omega\circ  \Lambda=\Omega\times_H \Lambda=\Omega\times \Lambda/\sim\]
where, for all  $x\in \Omega$,  $y\in \Lambda$ and $h\in H$, $(x,h y)\sim (x h,y)$.
\end{defi}

This construction is compatible with isomorphisms, and, endowed with the induced composition law, $\Bc$ defines a category. 

As we work with left-free bisets, we consider the subcategory $\A\subseteq \Bc$ where 
the objects are the same but we restrict the morphisms to isomorphism classes of left-free bisets.
This gives us a category and the composition follows from the next lemma.

\begin{lem}\label{comp}
Let $G,H$ and $K$ be finite groups.
Let $[K,\varphi]\in\A(G,H)$ and $[L,\psi]\in \A(H,K)$.
Then, \[[K,\varphi]\circ[L,\psi]=\coprod_{x\in \varphi(K)\backslash H/L} [\varphi^{-1}(\varphi(K)\cap  xLx^{-1}),\psi\circ c_{x^{-1}}\circ\varphi].\]
\end{lem}

\begin{proof}
We refer to that identity as the double coset formula and it is a direct consequence of \cite{Bo} Proposition 1. 
\end{proof}
\subsection{\texorpdfstring{$\Fc$}{Lg}-characteristic bisets and trivial coefficients}\label{4.2}

Let $(S,\Fc,\Li)$ be a $p$-local finite group.
When we work with trivial coefficients, the idea is to consider the category $\A_\Fc$ defined as follows. $\Ob(\A_\Fc)$ is the set of subgroups of $S$ and, for $P,Q\leq S$, $\A_\Fc(P,Q)$ is the set of isomorphism classes of $\Fc$-generated left-free $(P,Q)$-bisets, i.e. the $(P,Q)$-bisets union of transitive bisets of the form $[R,\varphi]$ with $R\leq P$ and $\varphi\in\Hom_\Fc(R,Q)$.

Then, for $M$ a $\Zp$-module, we construct a functor 
\[\xymatrix{M:\A_\Fc\ar[r] & \Zp\Mod}\]
defined on objects by $M(P)=H^*(P,M)$ for every $P\leq S$ and on morphisms as follows.
For every $P,Q\leq S$, $R\leq P$ and $\varphi\in\Hom_\Fc(P,Q)$, $({}_P[R,\varphi]_Q)_*=\tr_R^P\circ\varphi^*$. More generally, for every $\Fc$-generated left-free $(P,Q)$-biset $\Omega$ we define $\Omega_*$ by sum of its transitive components.

The existence of this functor will help us to construct an idempotent of $H^*(S,M)$ with image $H^*(\Fc^c,M)$.
For that, we also need the notion of $\Fc$-characteristic $(S,S)$-biset.

\begin{defi}\label{fcstable}
 Let $\Omega$ be a left-free $(S,S)$-biset.
\begin{enumerate}[(a)]
 \item We say that $\Omega$ is \emph{$\Fc$-generated} if it is the union of $(S,S)$-bisets of the form $[P,\varphi]$ with $P\in\Ob(\Fc)$ and $\varphi\in \Hom_\Fc(P,S)$.
 \item We say that $\Omega$ is \emph{left-$\Fc$-stable} if for all $P\in \Ob(\Fc)$ and $\varphi\in\Hom_\Fc(P,S)$, we have 
 $_\varphi\Omega_S\cong {}_P\Omega_S$, i.e.
 \[\left({}_P[P,\varphi]_S\right)\circ [\Omega]=\left({}_P[P,\incl^S_P]_S\right)\circ [\Omega].\] 
 \item We say that $\Omega$ is \emph{right-$\Fc$-stable} if for all $P\in \Ob(\Fc)$ and $\varphi\in\Hom_\Fc(P,S)$, we have 
 $_S\Omega_\varphi\cong {}_S\Omega_P$, i.e.
 \[[\Omega]\circ \left({}_S[\varphi(P),\varphi^{-1}]_P\right)  = [\Omega]\circ \left({}_S[P,\Id_P^S]_P\right).\]
 \item We say that $\Omega$ is \emph{non degenerate} if $|\Omega|/|S|\neq 0$ modulo $p$.
\end{enumerate}
If $\Omega$ satisfies all this four properties, we say that $\Omega$ is an $\Fc$-characteristic $(S,S)$-biset.
\end{defi}

The notion of $\Fc$-characteristic biset was first motivated by unpublished work of Linckelmann and Webb. 
They are the ones who first formulated these conditions and
recognized the importance of finding a biset with these properties. Broto, Levi and Oliver proved that such a biset always exists if the fusion system is saturated.

\begin{prop}[\cite{BLO2}, Proposition 5.5]
Let $\Fc$ be a fusion system over a $p$-group $S$.
If $\Fc$ is saturated, then there exists an $\Fc$-characteristic $(S,S)$-biset.
\end{prop}

In fact, Ragnarsson and Stancu (\cite{RS}, Theorem A), and independently Puig (\cite{P7}, Proposition 21.9), proved  that a fusion system $\Fc$ is saturated if, and only if, 
there exists a $\Fc$-characteristic $(S,S)$-biset.

Let $M$ be a $\Zp$-module, any $\Fc$-characteristic biset induces
an idempotent of $H^*(S,M)$ with image $H^*(\Fc^c,M)$.

\begin{prop}[cf. \cite{BLO2}, Proposition 5.5]\label{triv}
 Let $(S,\Fc,\Li)$ be a $p$-local finite group and $M$ be a $\Zp$-module (with a trivial action of $\pi_\Li$).
If $\Omega$ is an $\Fc$-characteristic biset, then $\frac{|S|}{|\Omega|}\Omega_*\in\End(H^*(S,M))$ defines an idempotent with image $H^*(\Fc^c,M)$.
\end{prop}
  
\begin{proof}
In \cite{BLO2}, Proposition 5.5, this is proved for $M=\F_p$ but the general case works the same way.
\end{proof}

\subsection{Bisets and twisted coefficients}\label{4.3}
Let $(S,\Fc,\Li)$ be a $p$-local finite group. When we work with twisted coefficients, one has to be more careful.
Unlike the case of trivial coefficients, defining a functor from $\A_\Fc$ to $\Zp\Mod$ does not work in general.

In fact, for $M$ a $\Zp[\pi_\Li]$-module, our cohomological functor $H^*(-,M)$ cannot be defined on $\Fc$ but only on $\Fc^c$ and thus, we can only consider $\Fc^c$-generated bisets. 

\begin{defi}
Let $P,Q$ be two $\Fc$-centric subgroups of $S$.
A left-free $(P,Q)$-biset is \emph{$\Fc^c$-generated} if it is an union of transitive bisets of the form $[R,\varphi]$ with $R\in\Ob(\Fc^c)$ and $\varphi\in\Hom_\Fc(R,Q)$. 
\end{defi}

Unfortunately, we can see from Lemma \ref{comp} that the set of isomorphism classes of $\Fc^c$-generated bisets is not stable with respect to composition.
Hence, we can not, by analogy with $\A_\Fc$, define a category $\A_{\Fc^c}$
where the objects are the $\Fc$-centric subgroups of $S$
and for $P$ and $Q$ two $\Fc$-centric subgroups of $S$, $\A_{\Fc^c}(P,Q)$ 
is the set of isomorphism classes of $\Fc^c$-generated left-free $(P,Q)$-bisets.


Nevertheless, for all $\Zp[\pi_\Li]$-module $M$ and $P,Q\leq S$, we have a map from the set $A_{\Fc^c}(P,Q)$ of isomorphism classes of $\Fc^c$-generated left-free $(P,Q)$-bisets to $\Hom(H^*(P,M),H^*(Q,M))$.

For $P,Q,R\in\Fc^c$ with $R\leq Q$ and $\varphi\in\Hom_\Fc(R,P)$,
we can associate to the $(P,Q)$-pair $\lbrace R,\varphi\rbrace$ a morphism 
\[\xymatrix@R=1mm{ \lbrace R,\varphi\rbrace_*=\tr_R^Q\circ\varphi^*: H^*(P,M)\ar[r] & H^*(Q,M).}\]

If we consider another $(P,Q)$-biset $\lbrace R',\varphi'\rbrace$ isomorphic to $\lbrace R,\varphi\rbrace$ 
(this implies that $R'$ is also $\Fc$-centric), we obtain the same morphism. 
%
Then we can set $[R,\varphi]_*$ as the composite $\tr_R^Q\circ\varphi^*$ and it is well-defined. 
Finally, for $\Omega$ a left-free $\Fc^c$-generated $(P,Q)$-biset, we define $\Omega_*$ by the sum of its transitive components.

\begin{rem}\label{Omegadeltafunctor}
 By Proposition \ref{phimorhdelta}, for $\varphi\in\Mor(\Fc^c)$, $\varphi^*=H^*(\varphi,-)$ is a morphism of $\delta$-functors.
 Hence, as $\Omega_*$ is a sum of composites of transfers, restrictions and $\varphi^*$, for $\varphi\in\Mor(\Fc^c)$,
 which are all morphisms of $\delta$-functors, it is an endomorphism of the $\delta$-functor $\left(H^*(S,-),\delta_{H^*(S,-)}\right)$.
\end{rem}

\subsection{Idempotents and twisted coefficients}\label{4.4}

In general, an $\Fc$-characteristic biset is not $\Fc^c$-generated.
Hence, when we are working with twisted coefficients, we cannot use directly an $\Fc$-characteristic biset as in the trivial case.
However, we can define, from an $\Fc$-generated $(S,S)$-biset, an endomorphism of $H^*(S,M)$ but not in a unique way. 
That is why we introduce the notion of $\Omega$-endomorphism.

\begin{defi}
 Let $(S,\Fc,\Li)$ be a $p$-local finite group.
 Let $P\leq S$ and $\varphi\in\Hom_\Fc(P,S)$. Let 
 \[
 \varphi= \incl_{\varphi(P)}^S\circ\pi(\psi_1)|^{\varphi(P)}_{P_{r-1}}\circ\dots\circ\pi(\psi_1)|^{P_1}_{P}
 \]
be a decomposition of $\varphi$ into automorphisms of $\Fc$-centric subgroups. 
 
Given a $\Zp[\pi_\Li]$-module $M$, a \emph{${}_S[P,\varphi]_S$-endomorphism} is an endomorphism of $H^*(S,M)$ given by the following composition
\begin{align*}
\tr_P^S \circ \left( H^*(\pi(\psi_1)|^{P_1}_{P},\omega(\psi_1)^{-1}) \right.&\circ  H^*(\pi(\psi_1)|^{P_2}_{P_1},\omega(\psi_2)^{-1})&\\
 &\circ\left. \cdots \circ H^*(\pi(\psi_1)|^{\varphi(P)}_{P_{r-1}},\omega(\psi_r)^{-1})\circ\Res_{\varphi(P)}^S \right).&
\end{align*}

 More generally, we define an \emph{$\Omega$-endomorphism}, for $\Omega$ an $\Fc$-generated $(S,S)$-biset as a 
 sum of the previous morphisms given by the transitive components of $\Omega$. 
 \end{defi}

\begin{rem}\label{omegatrivial}
If the action of $\pi_\Li$ on $M$ is trivial, every $\Omega$-endomorphism corresponds to $\Omega_*$ and, 
 if $\Omega$ is $\Fc$-characteristic, by Proposition \ref{triv}, $\frac{|S|}{|\Omega|}\Omega_*$ is an idempotent with image $H^*(\Fc^c,M)$.
\end{rem}

Let us look at the behavior of $\Omega$-endomorphisms induced by $\Fc$-characteristic $(S,S)$-bisets with $\Fc^c$-stable elements.

\begin{lem}\label{nontriv}
Let $(S,\Fc,\Li)$ be a $p$-local finite group and $M$ be a $\Zp[\pi_\Li]$-module.
 If $\Omega$ is an $\Fc$-characteristic $(S,S)$-biset and $\omega_*$ is an $\Omega$-endomorphism, 
 then $\frac{|S|}{|\Omega|}\omega_*\in\End(H^*(S,M))$ restricted to $H^*(\Fc^c,M)$ is the identity. 
\end{lem}

\begin{proof}
Let $P$ be a subgroup of $S$, $\varphi\in\Hom_\Fc(P,S)$ and $\lambda_*$ a ${}_S[P,\varphi]_S$-endomorphism.
For every $x\in H^*(\Fc^c,M)$, 
\[
\lambda_*(x)=\tr_P^S\circ \Res_P^S(x)=[S:P]x=\frac{|[P,\varphi]|}{|S|}x. 
\]
Hence, for every $x\in H^*(\Fc^c,M)$,
\[\omega_*(x)= \frac{|\Omega|}{|S|}x.\]
\end{proof}

\begin{rem}\label{omegadeltafunctor}
Notice also that for $\Omega$ an $\Fc$-characteristic $(S,S)$-biset,
 by construction and Remark \ref{phimorhdelta2}, an $\Omega$-endomorphism defines a morphism of $\delta$-functors.
\end{rem}

We remind the reader that in this article, every $p$-group is finite. This turn out to be a crucial property in the following result.

\begin{prop}\label{idem}
 Let $(S,\Fc,\Li)$ be a $p$-local finite group and let $M$ be an abelian $p$-group with an action of $\pi_\Li$. 
Let $\Omega$ be an $\Fc$-characteristic $(S,S)$-biset and $\omega_*$ an $\Omega$-endomorphism.
 For every $k\geq 0$, there is a natural number $N_{k,M}>0$ such that 
 $\left(\frac{|S|}{|\Omega|}\omega_k\right)^{N_{k,M}}$ defines an idempotent $\overline{\omega}_{k,M}$ of $H^k(S,M)$ and we have 
 \[H^k(\Fc^c,M)\subseteq \im(\overline{\omega}_{k,M}).\]
\end{prop}

\begin{proof}
To simplify the notation, we write $\omega=\frac{|S|}{|\Omega|}\omega_k$.
 For any $k\geq 0$, we have the following decreasing family of subgroups of $H^k(S,M)$.
\[
 H^k(\Fc^c,M)\subseteq \dots \subseteq \im(\omega^r)\subseteq\im(\omega^{r-1})\subseteq\dots \subseteq \im(\omega^1)\subseteq\im(\omega^0)=H^k(S,M).
\]
As $H^k(S,M)$ is a finite abelian $p$-group, this sequence stabilizes. 
Thus there is an $n_0\geq 1$ such that for all $n\geq n_0$ $\im(\omega^n)=\im(\omega^{n_0})$.
In particular, $\omega^{n_0}|_{\im(\omega^{n_0})}$ is a permutation of the finite set $\im(\omega^{n_0})$ and there is an $l$
such that $(\omega^{n_0}|_{\im(\omega^{n_0})})^l=\Id_{\im(\omega^{n_0})}$.
Thus, for $N_{k,M}=l\times n_0$, the endomorphism $\overline{\omega}_{k,M}=\omega^{N_{k,M}}\in\End(H^k(S,M))$ is an idempotent with image $\im(\omega^{n_0})\supseteq H^k(\Fc^c,M)$.
\end{proof}

Hence, we can define an idempotent of $H^*(S,M)$ as follows. For every $k\geq 0$ and every $x\in H^k(S,M)$, 
\[
\overline{\omega}_{k,M}(x)=\left(\frac{|S|}{|\Omega|}\omega_k\right)^{\prod_{i=0}^k N_{i,M}} (x).                                                                                                                    
 \]

Moreover, this definition only depends on the $\Omega$-endomorphism $\omega$.

\begin{defi}
For $\Omega$ an $\Fc$-characteristic $(S,S)$-biset and $\omega_*$ an $\Omega$-endomorphism,
the idempotent $\overline{\omega}_{*,-}$ of $H^*(S,-)$ obtained by the previous process is called the \emph{$\Fc^c$-characteristic idempotent} associated to $\omega$. 

Let $M$ be an abelian $p$-group with an action of $\pi_\Li$, we denote by $I_\omega^*(M)\subseteq H^*(S,M)$ the image of $\overline{\omega}_{*,M}$. 
\end{defi}

\begin{rem}\label{rk}
 Notice that, by Remark \ref{omegatrivial}, if the action on $M$ is trivial, then $I_\omega^*(M)=H^*(\Fc^c,M)$.
\end{rem}

\begin{prop}\label{Omegadelta}
 Let $(S,\Fc,\Li)$ be a $p$-local finite group. 
If $\Omega$ is an $\Fc$-characteristic $(S,S)$-biset and $\omega_*$ an $\Omega$-endomorphism,
 then $\overline{\omega}_{*,-}$, the $\Fc^c$-characteristic idempotent induced by $\omega_*$, 
 defines an endomorphism of the $\delta$-functor $\left(H^*(S,-),\delta_{H^*(S,-)}\right)$.
\end{prop}

\begin{proof}
For $M$ an abelian $p$-group with an action of $\pi_\Li$ and $k\geq 0$, we denote by $N_{k,M}$ a natural number as in Proposition \ref{idem}.

We have first to show that $\overline{\omega}_{*,-}$, the $\Fc^c$-characteristic idempotent associated to $\omega_*$,
 defines a natural transformation from the functor $H^*(S,-)$ to itself.
For every pair of abelian $p$-groups $(M,N)$ with an action of $\pi_\Li$ and every $\varphi\in\Hom_{\Zp[\pi_\Li]}(M,N)$, 
let us consider, for $k\geq 0$, the following diagram,
\[
\xymatrix{ 
 H^k(S,M) \ar[r]^\Id \ar[d]_{\overline{\omega}_{k,M}}  & H^k(S,M) \ar[r]^{\varphi_k} \ar[d]_{\tilde{\omega}_{k,M,N}} & H^k(S,N) \ar[d]_{\tilde{\omega}_{k,M,N}} \ar[r]^\Id & H^k(S,N) \ar[d]^{\overline{\omega}_{k,N}} \\
 H^k(S,M) \ar[r]_\Id         & H^k(S,M) \ar[r]_{\varphi_k}       & H^k(S,N)        \ar[r]_\Id & H^k(S,N)  
 }
\]
where $\tilde{\omega}_{k,M,N}=\left(\omega_k\right)^{\prod_{i=0}^k N_{i,M}\times \prod_{i=0}^k N_{i,N} }$.
The middle square commutes as $\tilde{\omega}_{k,M,N}$ is a finite iteration of $\frac{|S|}{|\Omega|}\omega_k$ 
and $\omega_*$ is an endomorphism of $\delta$-functors by Remark \ref{omegadeltafunctor}.
The leftmost square commutes because, as $\overline{\omega}_{k,M}$ is an idempotent of $H^k(S,M)$, 
$\tilde{\omega}_{k,M,N}=\overline{\omega}_{k,M}^{\prod_{i=0}^k N_{i,N}}=\overline{\omega}_{k,M}$. 
Finally, the rightmost one commutes because, as $\overline{\omega}_{k,N}$ is an idempotent of $H^k(S,N)$, 
$\tilde{\omega}_{k,M,N}=\overline{\omega}_{k,N}^{\prod_{i=0}^k N_{i,M}}=\overline{\omega}_{k,N}$.
Hence, the exterior diagram commutes.

Now, to show that it defines a morphism of $\delta$-functor, 
let us consider a short exact sequence of abelian $p$-groups with an action of $\pi_\Li$, $\xymatrix{0\ar[r] & L\ar[r] & M\ar[r] & N\ar[r] & 0}$.
By the previous argument we just have to show that, for $k\geq 0$, the following diagram commutes,
\[
 \xymatrix{
 H^k(S,N) \ar[r]^{\delta}\ar[d]_{\overline{\omega}_{k,N}} & H^{k+1}(S,L) \ar[d]^{\overline{\omega}_{k+1,L}} \\
 H^k(S,N) \ar[r]_\delta & H^{k+1}(S,L)
 }
\]
where $\delta=\delta_{H^*(S,-)}$ corresponds to the connecting homomorphism.
Consider then the following diagram,
\[
\xymatrix{ 
 H^k(S,N) \ar[r]^\Id \ar[d]_{\overline{\omega}_{k,N}}  & H^k(S,N) \ar[r]^{\delta} \ar[d]_{\tilde{\omega}_{k,L,N}} & H^{k+1}(S,L) \ar[d]_{\tilde{\omega}_{k+1,L,N}} \ar[r]^\Id & H^{k+1}(S,L) \ar[d]^{\overline{\omega}_{k+1,L}} \\
 H^k(S,N) \ar[r]_\Id         & H^k(S,N) \ar[r]_{\delta}       & H^{k+1}(S,L)        \ar[r]_\Id & H^{k+1}(S,L)  
 }
\]
where 
\[\tilde{\omega}_{k,L,N}=\left(\overline{\omega}_k\right)^{\prod_{i=0}^{k+1} N_{i,L}\times \prod_{i=0}^{k+1} N_{i,N} }\]
and
\[\tilde{\omega}_{k+1,L,N}=\left(\overline{\omega}_{k+1}\right)^{\prod_{i=0}^{k+1} N_{i,L}\times \prod_{i=0}^{k+1} N_{i,N} }.\]
The middle square commutes as $\tilde{\omega}_{k,L,N}$ and $\tilde{\omega}_{k+1,L,N}$ 
 are finite iterations of $\overline{\omega}_k$ and $\overline{\omega}_{k+1}$, and $\overline{\omega}_*$ 
 is an endomorphism of $\delta$-functors by Remark \ref{omegadeltafunctor}.
The leftmost square commutes because, as $\overline{\omega}_{k,N}$ is an idempotent of $H^k(S,N)$, 
$\tilde{\omega}_{k,L,N}=\overline{\omega}_{k,N}^{N_{k+1,N}\times\prod_{i=0}^{k+1} N_{i,L}}=\overline{\omega}_{k,N}$. 
The rightmost one commutes because, as $\overline{\omega}_{k,L}$ is an idempotent of $H^k(S,L)$, 
$\tilde{\omega}_{k+1,L,N}=\overline{\omega}_{k+1,L}^{\prod_{i=0}^{k+1} N_{i,N}}=\overline{\omega}_{k+1,L}$.
Thus, the exterior diagram commutes.
\end{proof}

\subsection{The idempotent for a constrained fusion system}\label{4.5}

When we work with a constrained fusion system, the $(S,S)$-characteristic biset is $\Fc^c$-generated and, working
with a suitable category, it induces, for every $\Zp[\pi_\Li]$-module $M$, an idempotent of $H^*(S,M)$ with image $H^*(\Fc^c,M)$.
Let us first recall the notion of constrained fusion system.

\begin{defi}\label{constrained}
Let $\Fc$ be a fusion system over a $p$-group $S$.

A subgroup $Q\leq S$ is \emph{normal in $\Fc$} if $Q\trianglelefteq S$, and for all $P,R\leq S$ and every $\varphi\in\Hom_\Fc(P,R)$, $\varphi$ extends to a morphism
$\overline{\varphi}\in\Hom_\Fc(PQ,RQ)$ such that $\overline{\varphi}(Q)=Q$. 

We write $O_p(\Fc)$ for the maximal subgroup of $S$ which is normal in $\Fc$.

We say that  $\Fc$ is \emph{constrained} if $O_p(\Fc)$ is $\Fc$-centric. 
\end{defi}

Define, for $\Fc$ a fusion system over a $p$-group $S$ and $P_0$ a subgroup of $S$, $\A_{\Fc\geq P_0}$
as follow.
 \[\Ob(\A_{\Fc\geq P_0})=\lbrace P_0\leq P\leq S\rbrace\text{ is the set of all subgroups of $S$ containing $P_0$}\]
and for all $P,Q\in\Ob(\A_{\Fc\geq P_0})$,
 \[\A_{\Fc\geq P_0}(P,Q)=\lbrace\text{$\Fc$-generated left-free $(P,Q)$-bisets union of $[R,\varphi]$ with $R\geq P_0$}\rbrace.\]

$A_{\Fc\geq P_0}$ is not in general a subcategory of $\A_\Fc$. 
The problem comes from Lemma \ref{comp}: 
the set \[\Mor(\A_{\Fc\geq P_0})=\bigsqcup_{P,Q\in\Ob(\A_{\Fc\geq P_0})}\A_{\Fc\geq P_0}(P,Q)\] is not stable with respect to composition.
But it is stable when the subgroup $P_0\leq S$ is \textit{weakly closed} in $\Fc$, i.e. 
$P^\Fc=\{P\}$.
%

\begin{lem}
 Let $\Fc$ be a fusion system over a $p$-group $S$.
 If $P_0\trianglelefteq S$ is weakly closed in $\Fc$, then $\A_{\Fc\geq P_0}$, with the composition defined in \ref{compdef}, is a subcategory of $\A_\Fc$.
\end{lem}

\begin{proof}
 As $P_0$ is weakly closed in $\Fc$, for every $R,P\geq P_0$, $s\in S$ and $\varphi\in\Hom_{\Fc}(R,S)$, 
\[\varphi^{-1}(\varphi(R)\cap sPs^{-1})\geq \varphi^{-1}(\varphi(P_0)\cap sP_0s^{-1})=P_0.\]  
Thus, by Lemma \ref{comp}, $\Mor(\A_{\Fc\geq P_0})$ is stable with respect to composition
 and $\A_{\Fc\geq P_0}$ defines a subcategory of $\A_\Fc$.
\end{proof}

For example, the subgroup $O_p(\Fc)$ is normal in $\Fc$. Thus it is weakly closed in $\Fc$ and $\A_{\Fc\geq O_p(\Fc)}$ is a subcategory of $\A_\Fc$.
\vspace{10pt}

When $\Fc$ is constrained, $O_p(\Fc)$ is $\Fc$-centric. Thus, every biset $\Omega\in\Mor(\A_{\Fc\geq O_p(\Fc)})$ is $\Fc^c$-generated.
Hence, if $\Fc$ is constrained, for every $\Zp[\pi_\Li]$-module $M$, we have, as in the trivial case, a functor 

\[
 \xymatrix@R=1mm{\A_{\Fc\geq O_p(\Fc)}\ar[r] & \Zp\Mod \\
                     P\ar@{|->}[r] & H^*(P,M) \\
                     {}_P[R,\varphi]_Q\ar@{|->}[r]& \tr_R^P\circ\varphi^*.}
\]

Moreover, if we look at the \emph{minimal} $\Fc$-characteristic $(S,S)$-biset i.e. the smallest $\Fc$-characteristic $(S,S)$-biset, we have the following.

\begin{prop}\label{constrained biset}
 Let $\Fc$ be a constrained fusion system over a $p$-group $S$.
 If $\Omega$ is the minimal $\Fc$-characteristic biset, then $\Omega\in\A_{\Fc\geq O_p(\Fc)}$.
\end{prop}

\begin{proof}
 This is a direct corollary of \cite{GRh}, Proposition 9.11.
 Indeed, by \cite{GRh}, Proposition 9.11, every $[P,\varphi]$ which appears in the decomposition of $\Omega$ satisfies $P\geq O_p(\Fc)$.
\end{proof}

Hence, using the same argument as for Proposition \ref{triv}, we have the following theorem.

\begin{thm}
 Let $(S,\Fc,\Li)$ be a $p$-local finite group and $M$ be a $\Zp[\pi_\Li]$-module.
 Let $\Omega$ be the minimal $\Fc$-characteristic $(S,S)$-biset. 
 If $\Fc$ is a constrained fusion system, then $\frac{|S|}{|\Omega|}\Omega_*\in\End(H^*(S,M))$ is an idempotent with image the $\Fc^c$-stable elements $H^*(\Fc^c,M)$.
\end{thm}

\begin{proof}
 By Proposition \ref{constrained biset}, $\Omega\in\A_{\Fc\geq O_p(\Fc)}$, and the proof is the same as the proof of Proposition \ref{triv}.
\end{proof}

\subsection{A \texorpdfstring{$\delta$}{Lg}-functor}\label{4.6}

Let $(S,\Fc,\Li)$ be a $p$-local finite group, $\Omega$ be an $\Fc$-characteristic $(S,S)$-biset,
and $\omega_*$ an $\Omega$-endomorphism.

For $M$ an abelian $p$-group with an action of $\pi_\Li$,
let $\overline{\omega}_{*,-}\in\End(H^*(S,M))$ be the associated $\Fc^c$-characteristic idempotent.

Let us start with the behavior of $\delta$-functors with idempotents. 
We recall that a $\delta$-functor can be seen as a functor from the category $\Ss_\A$ of short exact sequences in $\A$ to $\Ch(\B)$, the category of $\mathbb{Z}$-graded chain complexes in $B$, which sends any short exact sequence to an acyclic chain complex. We refer the reader to \cite{We} for more details and properties. 

\begin{lem}\label{imexact}
 Let 
 $(M_*,f_*)=\left(\cdots\xrightarrow{f_{l-2}}M_{l-1}\xrightarrow{f_{l-1}} M_l\xrightarrow{f_l} M_{l+1}\xrightarrow{f_{l+1}} \cdots\right)_{l\in\mathbb{Z}}$ be
a long exact sequence in an abelian category $\A$.
Let $i_*:(M_*,f_*)\rightarrow (M_*,f_*)$ be a morphism of long exact sequences such that, for all $l\in\mathbb{Z}$, $i_l$ is an idempotent of $M_l$.
Then the sequence \[\cdots\xrightarrow{f_{l-2}}\im(i_{l-1})\xrightarrow{f_{l-1}}\im(i_l)\xrightarrow{f_l} \im(i_{l+1})\xrightarrow{f_{l+1}} \cdots\]
is exact.
\end{lem}

\begin{proof}
Let $l\in \mathbb{Z}$ and $x\in \im(i_l)$ such that $f_l(x)=0$. By exactness of $(M_*,f_*)$ in $l$, 
there is a $y\in M_ {l-1}$ such that $f_{l-1}(y)=x$.
Thus $x=i_l(x)=i_l\circ f_{l-1}(y)=f_{l-1}\circ i_{l-1}(y)$ and hence we obtain the exactness of $(\im(i_*),f_*)$ in degree $l$.
\end{proof}

\begin{prop}\label{deltaidem}
Let $\A,\B$ be two abelian categories and let $\left(F^*,\delta_F\right):\A\rightarrow\B$ be a $\delta$-functor.
If $\xymatrix{i^*:\left(F^*,\delta_F\right)\rightarrow \left(F^*,\delta_F\right)}$ 
is an idempotent of $\delta$-functors, then $\left(\im(i^*),\delta_F\right)$ defines a $\delta$-functor.    
\end{prop}

\begin{proof}
 A $\delta$-functor can be seen as a functor from the category $\Ss_\A$ of short exact sequences in $\A$ to $\Ch(\B)$ 
 which sends any short exact sequence to an acyclic chain complex. 
 A morphism of $\delta$-functors is then a natural transformation in that setting.
 With this point of view, this is just a corollary of Lemma \ref{imexact}.
\end{proof}

\begin{thm}\label{deltafunctor}
 Let $(S,\Fc,\Li)$ be a $p$-local finite group, $\Omega$ an $\Fc$-characteristic $(S,S)$-biset
 and $\omega_*$ an $\Omega$-endomorphism.
 Then, the functor $I_\omega^*(-)$, with the connecting homomorphism $\delta_{H^*(S,-)}$,
 defines a $\delta$-functor from the category of finite $\Zp[\pi_\Li]$-modules to $\Zp\Mod$.
\end{thm}

\begin{proof}
 This is a direct corollary of Proposition \ref{Omegadelta} and Proposition \ref{deltaidem}.
\end{proof}

In the next section, we will show that if the action on $M$ is nilpotent, then $I_\omega^*(M)=H^*(\Fc^c,M)$. 
But this is not clear at all in general. This raised the following question.

\begin{conj}
Let $(S,\Fc,\Li)$ be a $p$-local finite group, $\Omega$ an $\Fc$-characteristic $(S,S)$-biset and $\omega$ an $\Omega$-endomorphism. If $M$ is an abelian $p$-group with an action of $\pi_1(\pcompl{|\Li|})$, then
\[H^*(\Fc^c,M)\cong I_\omega^*(M).\]
\end{conj}

We insist that, in view of the counterexamples given by Levi and Ragnarsson (\cite{LR} Proposition 3.1), we cannot expect $I_\omega^*(M)$ to be isomorphic to $H^*(|\Li|,M)$ in general. But it is natural to ask if $I_\omega(M)$ always corresponds to the $\Fc^c$-stable elements. 

\section{The cohomology of the geometric realization of a linking system with nilpotent coefficients}\label{5}

We give here a proof of the main theorem.

\begin{lem}\label{lem2}
 Let $(S,\Fc,\Li)$ be a $p$-local finite group.
 The natural inclusion $\delta_S$ of $\Bc(S)$ in $\Li$ induces, for any $\Zp[\pi_\Li]$-module $M$, a natural morphism in cohomology
             \[ \xymatrix{ H^*(|\Li|,M)\ar[r] & H^*(\Fc^c,M)\subseteq H^*(S,M).}\]
\end{lem}

\begin{proof}
This follows easily from the functoriality of the geometric realization.
\end{proof}

\begin{lem}\label{lem}
 Let $(S,\Fc,\Li)$ be a $p$-local finite group and let $\Omega$ be an $\Fc$-characteristic $(S,S)$-biset.
Let also $0\rightarrow L \rightarrow M \rightarrow N \rightarrow 0$ be a short exact sequence of finite $\Zp[\pi_\Li]$-modules.
If $\delta_S$ induces isomorphisms $ H^*(|\Li|,L)\cong I_\omega^*(L)$ and $ H^*(|\Li|,N)\cong I^*_\omega(N)$,
then $\delta_S$ induces an isomorphism \[ H^*(|\Li|,M)\cong I_\omega^*(M).\]  
\end{lem}

\begin{proof}
 Consider the exact sequences in cohomology induced by the short exact sequence
\[\xymatrix{0 \ar[r] & L \ar[r] & M\ar[r] & N \ar[r] & 0}\]
and look at the following diagram (where $\overline{\omega}_{*,-}$ denote the $\Fc^c$-characteristic idempotent associated to $\Omega$).
\[
\xymatrix{\cdots\ar[r] & H^{n-1}(|\Li|,N)\ar[r]\ar[d]_{\overline{\omega}_{n-1,N}\circ\delta_S^*} & H^{n}(|\Li|,L)\ar[r]\ar[d]_{\overline{\omega}_{n,L}\circ\delta_S^*} & H^{n}(|\Li|,M)\ar[r]\ar[d]_{\overline{\omega}_{n,M}\circ\delta_S^*} & H^{n}(|\Li|,N)\ar[r]\ar[d]_{\overline{\omega}_{n,N}\circ\delta_S^*} &  \cdots \\
\cdots\ar[r] & I^{n-1}_\Omega (N)\ar[r] & I^{n}_\Omega (L)\ar[r] & I^{n}_\Omega (M)\ar[r] & I^{n}_\Omega (N)\ar[r] & \cdots}
\]
As $H^*(|\Li|,-)$ is a $\delta$-functor and, by Theorem \ref{deltafunctor}, $I^{*}_\Omega$ is also a $\delta$-functor,
the two lines are exact and, as by Proposition \ref{Omegadelta} the $\Fc$-characteristic idempotent associated to $\Omega$
defines a morphism of $\delta$-functors, this diagram is commutative.
An application of the Five Lemma then finishes the proof.
\end{proof}

We can now state the main theorem.
\begin{thm}\label{main}
 Let $(S,\Fc,\Li)$ be a $p$-local finite group.
 If $M$ is an abelian $p$-group with a nilpotent action of $\pi_1(|\Li|)$, then $\delta_S$ induces a natural isomorphism
\[ H^*(|\Li|,M)\cong H^*(\Fc^c,M).\]
\end{thm}

\begin{proof}
 As the action of $\pi_\Li$ is nilpotent, there is a sequence 
 \[0=M_0\subseteq M_1\subseteq\dots\subseteq M_n=M\]
 such that, for every $1\leq i \leq n$, the action of $\pi_\Li$ on $M_{i}/M_{i-1}$ is trivial.
 We know, by Theorem \ref{trivial} and Remark \ref{rk}, that for $1\leq i \leq n$,  $\delta_S$ induces an isomorphism 
\[H^*(|\Li|,M_{i}/M_{i-1})\cong H^*(\Fc^c,M_{i}/M_{i-1})=I^*_\Omega(M_{i}/M_{i-1}).\]
 By induction on $n$, and by Lemma \ref{lem}, we get that $H^*(|\Li|,M)\cong I_\omega^*(M)$.
 Finally we also have by Lemma \ref{lem2} that \[\delta_S(H^*(|\Li|,M))\subseteq H^*(\Fc^c,M)\subseteq I_\omega^*(M).\]
 Then $H^*(|\Li|,M)\cong H^*(\Fc^c,M)= I^*_\Omega(M)$.
\end{proof}

\section{The cohomology with twisted coefficients of \texorpdfstring{$p$}{Lg}-good spaces}\label{6}
We finish with a result on the cohomology with twisted coefficients of the $p$-completion of a $p$-good space
and we apply it, with Theorem \ref{main}, to compare the 
cohomology with twisted coefficients of $|\Li|$ and the $\Fc^c$-stable elements.

We refer the reader to Bousfield and Kan \cite{BK} for more details about $p$-completion. 
There is also a brief introduction in Aschbacher, Kessar and Oliver \cite{AKO}.
Here, for $X$ a space, 
\[\xymatrix{\lambda_X:X\rightarrow \pcompl{X}}\]
denote the structural natural transformation of the $p$-completion and we recall that if $X$ is $p$-good, it induces an isomorphism
 \[
 H^*(\pcompl{X},\F_p)\cong H^*(X,\F_p).
\]

\begin{lem}\label{5lemmapcomp}
 Let $X$ be a space and let $0\rightarrow L \rightarrow M \rightarrow N \rightarrow 0$ be a short exact sequence of $\Zp[\pi_1(\pcompl{X})]$-modules.
If $\lambda_X$ induces isomorphisms $ H^*(\pcompl{X},L)\cong H^*({X},L)$ and $ H^*(\pcompl{X},N)\cong H^*({X},N)$,
then $\lambda_X$ induces an isomorphism \[ H^*(\pcompl{X},M)\cong H^*({X},M).\]  
\end{lem}

\begin{proof}
This is a straightforward application of the Five Lemma.
\end{proof}

\begin{prop}\label{pgoodspace}
Let $X$ be a space and $M$ be an abelian $p$-group with an action of $\pi_1(\pcompl{X})$.
If $X$ is $p$-good and $\pi_1(\pcompl{X})$ is a finite $p$-group, then $\lambda_X$ induces a natural isomorphism
\[
 H^*(\pcompl{X},M)\cong H^*(X,M).
\]
\end{prop}

\begin{proof}
As $X$ is $p$-good, $\lambda_X$ induces an isomorphism, 
$H^*(\pcompl{X},\F_p)\cong H^*(X,\F_p).$
Moreover, as $\pi_1(\pcompl{X})$ is a $p$-group quotient of $\pi_1(X)$, the action of $\pi_1(\pcompl{X})$ on $M$ is nilpotent: 
there is a sequence 
 \[\{0\}=M_0\subseteq M_1\subseteq\dots\subseteq M_n=M\]
such that, for any $1\leq i \leq n$, $M_{i}/M_{i-1}\cong \F_p$ is the trivial module.
We conclude by induction on $n$ using Lemma \ref{5lemmapcomp}. 
\end{proof}

\begin{cor}\label{pgood}
Let $(S,\Fc,\Li)$ be a $p$-local finite group.
If $M$ is an abelian $p$-group with an action of $\pi_1(\pcompl{|\Li|})$, $\lambda_{|\Li|}$ induces an isomorphism
\[H^*(\pcompl{|\Li|},M)\cong H^*(|\Li|,M).\]
 
\end{cor}

\begin{proof}
 As $|\Li|$ is a $p$-good space and $\pi_1(\pcompl{|\Li|})$ is a finite $p$-group (\cite{AKO}, Theorem III.4.17), we can apply Proposition \ref{pgoodspace}.
\end{proof}

\begin{cor}\label{main2}
Let $(S,\Fc,\Li)$ be a $p$-local finite group.
If $M$ is an abelian $p$-group with an action of $\pi_1(\pcompl{|\Li|})$, then $\lambda_{|\Li|}\circ\delta_S^*$ induces a natural isomorphism
\[
 H^*(\pcompl{|\Li|},M)\cong H^*(\Fc^c,M).
\]
\end{cor}

\begin{proof}
By \cite{AKO}, Theorem III.4.17, $|\Li|$ is $p$-good and  $\pi_1(\pcompl{|\Li|})$ is a $p$-group. In particular, the action of $\pi_1(\pcompl{|\Li|})$  on $M$ is nilpotent.
Hence, this is just a corollary of Theorem \ref{main} and Corollary \ref{pgood}.
\end{proof}

%
%
%
\bibliography{biblio}{}
\bibliographystyle{abstract}



\end{document}